\documentclass{article}
\usepackage[style=alphabetic, backend=bibtex, maxnames=99]{biblatex}
\usepackage{amsfonts}
\usepackage{amsmath}
\usepackage{amssymb}
\usepackage{mathrsfs}
\usepackage{amsthm}
\usepackage{tikz}
\usetikzlibrary{arrows,positioning,decorations.pathmorphing,  decorations.markings}
\usepackage{tikz-cd}
\usepackage{indentfirst}
\usepackage{enumitem}
\usepackage{comment}
\usepackage{titlesec}
\usepackage{tocloft}
\usepackage{float}
\usepackage[linktocpage]{hyperref}
\usepackage[a4paper, total={6in, 8in}]
{geometry}

\newcommand{\R}{\mathbb{R}}
\newcommand{\C}{\mathbb{C}}
\newcommand{\N}{\mathbb{N}}
\newcommand{\Z}{\mathbb{Z}}
\newcommand{\Q}{\mathbb{Q}}

\newcommand{\A}{\mathbb{A}}

\newcommand{\T}{\mathbb{T}}

\newcommand{\bbS}{\mathbb{S}}
\newcommand{\bbk}{\Bbbk}

\newcommand{\cA}{\mathcal{A}}
\newcommand{\cB}{\mathcal{B}}
\newcommand{\cC}{\mathcal{C}}
\newcommand{\cE}{\mathcal{E}}

\newcommand{\cI}{\mathcal{I}}
\newcommand{\cK}{\mathcal{K}}
\newcommand{\cL}{\mathcal{L}}

\newcommand{\cU}{\mathcal{U}}
\newcommand{\cO}{\mathcal{O}}
\newcommand{\cP}{\mathcal{P}}
\newcommand{\cW}{\mathcal{W}}
\newcommand{\cX}{\mathcal{X}}

\newcommand{\fb}{\mathfrak{b}}
\newcommand{\fg}{\mathfrak{g}}

\newcommand{\fD}{\mathfrak{D}}
\newcommand{\fT}{\mathfrak{T}}

\newcommand{\sL}{\mathscr{L}}

\newcommand{\sP}{\mathscr{P}}

\newcommand{\Hom}{\textnormal{Hom}}

\newcommand{\Spec}{\textnormal{Spec }}

\newcommand{\x}{\times}
\newcommand{\ox}{\otimes}

\newcommand{\tr}{\textnormal{tr}}

\newcommand{\rank}{\textnormal{rank }}

\newcommand{\wt}{\widetilde}
\newcommand{\ve}{\varepsilon}

\newcommand{\uf}{\textnormal{uf}}
\newcommand{\s}{\textnormal{\textbf{s}}}
\newcommand{\rat}{\dashrightarrow}

\newcommand{\bfc}{\textnormal{\textbf{c}}}
\newcommand{\bfe}{\textnormal{\textbf{e}}}
\newcommand{\bff}{\textnormal{\textbf{f}}}
\newcommand{\bfg}{\textnormal{\textbf{g}}}
\newcommand{\bfi}{\textnormal{\textbf{i}}}

\newcommand{\bfv}{\textnormal{\textbf{v}}}

\newcommand{\bfC}{\textnormal{\textbf{C}}}
\newcommand{\bfD}{\textnormal{\textbf{D}}}
\newcommand{\bfE}{\textnormal{\textbf{E}}}
\newcommand{\bfF}{\textnormal{\textbf{F}}}

\newcommand{\bfK}{\textnormal{\textbf{K}}}

\newcommand{\diag}{\textnormal{diag}}
\newcommand{\prin}{\textnormal{prin}}
\newcommand{\Ad}{\textnormal{Ad}}
\newcommand{\Frac}{\textnormal{Frac}}

\newcommand{\trop}{\textnormal{trop}}

\newtheorem{thm}{Theorem}[section]
\newtheorem*{mainthm}{Main Theorem}
\newtheorem{lem}[thm]{Lemma}
\newtheorem{prop}[thm]{Proposition}

\newtheorem{cor}[thm]{Corollary}

\theoremstyle{definition}
\newtheorem{defn}[thm]{Definition}
\newtheorem{eg}[thm]{Example}
\newtheorem{remark}[thm]{Remark}

\newcommand{\Eq}[1]{\begin{align}#1\end{align}}

\tikzset{>=latex}
\tikzstyle{vthick}=[line width=1.8pt]

\newcommand\drawpath[2]{%
  \foreach \too [count=\c from 1] in {#1}
  {
  \ifthenelse{\c=1}
  {\xdef\from{\too}}
  {\path (\from) edge [->, #2] (\too);
    \xdef\from{\too}}
  };
}

\begin{document}

\title{On the polynomiality conjecture of cluster realization of quantum groups}

\author{  Ivan Chi-Ho Ip\footnote{
	  Department of Mathematics, Hong Kong University of Science and Technology\newline
	  Email: ivan.ip@ust.hk\newline		
	  Email: jyye@connect.ust.hk\newline
The first author is supported by the Hong Kong RGC General Research Funds GRF \#16305122.
          }, Jeff York Ye\footnotemark[1]
}
\maketitle

\numberwithin{equation}{section}
\begin{abstract}
In this paper, we give a sufficient and necessary condition for a regular element of a quantum cluster algebra $\cO_q(\cX)$ to be universally polynomial. This resolves several conjectures by the first author on the polynomiality of the cluster realization of quantum group generators in different families of positive representations.
\end{abstract}

\tableofcontents

\addtocontents{toc}{\protect\setcounter{tocdepth}{2}}

\section{Introduction}
\label{sec:intro}

Let $\fg$ be a semisimple Lie algebra. The quantum group $\cU_q(\fg)$ was originally defined by Drinfeld \cite{Drinfeld} and Jimbo \cite{Jimbo}, which is a Hopf algebra deforming the universal enveloping algebra $\cU(\fg)$. The Drinfeld's double $\fD(\fg)$ corresponds to the double of the Borel part $\cU_q(\fb)$. To study the structure of quantum groups, one can embed them into simpler algebras such as a quantum torus algebra $\C_q[\T]$ or its skew field of rational functions $\C_q(\T)$. This was explored by Gerasimov--Kharchev--Lebedev--Oblezin in \cite{GKL05} and \cite{GKLO} by generalizing the Gelfand--Tsetlin representations originally defined for classical Lie groups, and they obtained an embedding of $\cU_q(\mathfrak{sl}_n)$ into $\C_q(\T)$. The Feigin's homomorphism also provides an embedding of the positive Borel part $\cU_q(\fb_+)$ to the quantum torus algebra \cite{Ber96}. However, the embedding of the whole quantum group $\cU_q(\fg)$ into $\C_q[\T]$ is more difficult. For type $A_n$, this was solved by Kashaev--Volkov \cite{KV98}.\\

Meanwhile, these embeddings are related to quantum Liouville theory in physics, where the Chevalley generators of $\cU_q(\fg)$ act as positive self-adjoint operators on $L^2(\R)$. The $\cU_q(\mathfrak{sl}_2)$ case was studied by Schm\"udgen \cite{Sch96}, Ponsot--Teschner \cite{PT99}\cite{PT01}.\\

To study this problem for general Lie types, Frenkel and the first author introduced the notion of positive representations for split real quantum groups in \cite{FI14}. Split real quantum groups were developed by Faddeev in \cite{Fad95} and extended to its modular double in \cite{Fad99}, motivated from string theory and conformal field theory. Based on the idea by Schrader--Shapiro in the case of type $A_n$ \cite{SS19}, by forgetting the real structure the positive representations lead to an embedding of the quantum group $\cU_q(\fg)$ of general Lie type into a quotient of a certain quantum torus algebra $\cO_q(\cX_{\bfD(\bfi_0)})$ associated to a quiver $\bfD(\bfi_0)$ \cite{Ip18}. The quiver $\bfD(\bfi_0)$ is mutation equivalent to the quiver for the moduli space of framed $G$-local systems $\sP_{G,\odot}$, where $\odot$ denotes the once punctured disk with two special points on the boundary \cite{GS}. In later works, this embedding was modified in \cite{Ip20} and \cite{IM22} to give positive representations which are related to the quantization of the classical parabolic induction, and certain degeneration of the generalized Casimir operators.\\

It was conjectured in \cite[Section 1.1]{Ip18} and \cite[Conjecture 7.1]{Ip20} that the images of Chevalley generators of $\cU_q(\fg)$ are polynomials in the cluster variables for any cluster in the mutation class of the quantum torus algebra. It is also suspected that the same statement is true for the homomorphisms introduced in \cite{IM22}. The main goal of this paper is to resolve these conjectures.

\begin{mainthm}[Theorem \ref{MainThm}]
\label{MainThmInIntro}
For the homomorphisms
\begin{align*}
\kappa&: \fD_q(\fg) \to \cO_q(\cX)\\
\overline{\kappa}&: \fD_q(\fg) \to \cO_q(\overline{\cX})\\
\kappa_0&:\fD_q(\fg) \to \cO_q(\cX^0)
\end{align*}
described in Section \ref{sec:repn:cluster}, the image of the Chevalley generators $\bfE_i,\bfF_i,\bfK_i,\bfK_i'$ are all polynomials in the quantum cluster variables with coefficients in $\N[q^{\pm1}]$ for any cluster in the mutation class.
\end{mainthm}

The idea of the proof is to use $c$ and $g$-vectors, which are known to be sign coherent for skew-symmetric seed data \cite{GHKK}. This allows us to control the degree of the cluster $\cX$-variables after each mutation, and leads directly to the proof.\\

It will be important to understand the geometric meaning of the polynomiality of the image of quantum groups in the quantum cluster algebra. Just as the ring of Laurent polynomials $\C[X_k,X_k^{-1}]$ describes regular functions on the \emph{torus}, the polynomial ring $\C[X_k]$ is usually associated to regular functions on the \emph{affine space}. Therefore we can expect an interpretation of the quantum groups to be generalized from elements of the quantum torus algebra to elements associated to \emph{quantum plane}, motivated by \cite{Ip13}, in which $C^*$-and von Neumann algebraic techniques can be utilized to study the multiplier Hopf structures of their images, thus building a new connection between Drinfeld--Jimbo type quantum groups and non-commutative geometry. \\

Furthermore, it will be interesting to understand the combinatorial description of the polynomial images, where the coefficients are expected to count certain canonical objects associated to a given cluster seed. In \cite{SS17}, this is answered explicitly for type $A_n$ quantum groups in the case of the standard embedding, in which the coefficients of the polynomial images of the Borel part $\cU_q(\fb)$ gives a weighted count of the homology cycles around the dual planar bipolar graph of the associated basic quivers. The Lusztig's braid group action $T_{\bfi_0}$ described in \cite{GS} can be interpreted as a sequence of mutations which provides an $S_3$ action to the quantum group embedding for each triangle on the surface, and the results of this paper shows that these embeddings remain a polynomial. For example, in type $E_8$ and starting with the Coxeter words as initial seed, by a heavy brute force check, the $\bff_0$ generator is mutated up to a polynomial of 825887337 terms before reducing to the final configuration, with very curious coefficients in $\N[q^{\pm1}]$ represented in terms of ratios of $q$-numbers \cite{Ip20}. It will be important to give a more combinatorial description of these polynomials and their coefficients.\\

Finally, in \cite{IM22}, we proposed a program to classify the \emph{regular positive representations} as described in Section \ref{sec:discussion}. One of the defining conditions is that we have a homomorphism of the quantum group into universally Laurent polynomials of some quantum cluster variety. The results of this paper shows that all the classes considered so far have the image of quantum group lies in the \emph{universally polynomial} part of the quantum cluster variety, and it will be essential to understand whether this is characterization of such class of representations, or there are examples in between that is yet to be discovered.\\

The paper is organized as follows. In Section \ref{sec:prelim}, we give the definition of cluster varieties and quantum cluster algebra in general. In Section \ref{sec:repn}, we give the definition of quantum groups. Then we provide a summary of the works by the first author and Man \cite{Ip18}, \cite{Ip20}, \cite{IM22}, as well as the works by Goncharov and Shen \cite{GS}, \cite{She22} on embeddings of $\cU_q(\fg)$ into cluster algebras. In Section \ref{sec:polynomial}, we use the sign coherence of $c$-vectors to give a sufficient condition for a regular element of a quantum cluster algebra $\cO_q(\cX)$ to be universally polynomial. Using this, we prove the polynomiality conjectures by the first author. Another more abstract, but sufficient and necessary condition is also formulated using $G$-matrices. For the polynomiality of zero Casimir representations, we define cluster folding via a group action and relate the quantum cluster algebras $\cO_q(\cX)$ and $\cO_q(\cX^G)$. Finally, we provide some generalizations of the main theorem, including the amalgamation of multiple copies of the quiver and corresponds to the coproduct in $\cU_q(\fg)$. In Section \ref{sec:discussion}, we give some comment on the notion of regular positive representations introduced by Ip and Man in \cite{IM22}, in view of the new results obtained in Section  \ref{sec:polynomial}.

\section{Preliminaries on cluster algebras}
\label{sec:prelim}

\subsection{Basic definitions and properties}
\label{sec:prelim:defn}

We mainly follow the notation in \cite{FG2} and \cite{GHK13} for the $\cA$ and $\cX$ cluster varieties. However, we slightly modify the definitions following \cite[Section 6]{Nak21}, so that the mutation formulas for the $A$ and $X$ variables match that used by the first author in \cite{Ip18}, \cite{Ip20} and also the original definitions in \cite{FZ02}. (See Remark \ref{NotationIssue}.)

\begin{defn}\label{fixeddata}
A \textit{fixed data} $\Gamma$ consists of 
\begin{itemize}
\item A finite rank lattice $N$ with a skew-symmetric bilinear form
\Eq{
\{\cdot,\cdot\}:N\x N\to\Q.
}
\item An \textit{unfrozen} sublattice $N_\uf\subset N$ which is \textit{saturated}, i.e. if $kn\in N_\uf$ for some $k\in\Z$ and $n\in N$, then $n\in N_\uf$. If $N_\uf=N$, we say the fixed data has no frozen variables.
\item An index set $I$ with $|I|=n:=\rank N$ and a subset $I_\uf$ with $|I_\uf|=\rank N_\uf$.
\item Positive integers $d_i$ with $i\in I$ with greatest common divisor 1.
\item A sublattice $N^\circ\subset N$ of finite index such that $\{N^\circ,N_\uf\}\subset \Z$ and $\{N_\uf\cap N^\circ,N\}\subset \Z$. This sublattice $N^\circ$ depends on $d_i$ as described in Definition \ref{defseed} below.
\item Lattices $M:=\Hom(N,\Z)$ and $M^\circ :=\Hom(N^\circ,\Z)$.
\end{itemize}
\end{defn}

\begin{defn}\label{defseed}
For any fixed data $\Gamma$, a \textit{seed} $\s:=(\bfe_i)_{i\in I}$ is a labeled collection of elements of $N$ such that $\{\bfe_i\}_{i\in I}$ is a $\Z$-basis of $N$, $\{\bfe_i\}_{i\in I_\uf}$ is a $\Z$-basis of $N_\uf$, and $\{d_i\bfe_i\}_{i\in I}$ is a $\Z$-basis of $N^\circ$.
\end{defn}

Given a seed $\s$, we obtain a dual basis $\{\bfe_i^*\}$ for $M$, and a basis $\{\bff_i\}$ of $M^\circ$ given by \Eq{\bff_i:=d_i^{-1}\bfe_i^*.} We define the exchange matrix by \Eq{\ve_{ij}:=d_i\{ \bfe_i,\bfe_j\},} which belongs to $\Z$ if at least one of $i,j$ belong to $I_\uf$. The exchange matrix is skew-symmetrizable with left skew-symmetrizer $\diag(d_i^{-1})$.

The exchange matrix corresponds to a quiver $Q$ whose vertices are labeled by $I$. For each $i,j\in I$, if $\{\bfe_i,\bfe_j\}>0$, there are $\{\bfe_i,\bfe_j\}$ arrows from $i$ to $j$, and if $\{\bfe_i,\bfe_j\}<0$, there are $-\{\bfe_i,\bfe_j\}$ arrows from $j$ to $i$.

\begin{defn}
The map \Eq{n\mapsto \{\cdot,n\}} is well defined as maps $p_1: N_\uf\to M^\circ$ and $p_2: N\to M^\circ/N_\uf^\perp$ from the conditions $\{N^\circ,N_\uf\}\subseteq \Z$ and  $\{N_\uf\cap N^\circ,N\}\subseteq \Z$ respectively. Here taking dual of $N_\uf\cap N^\circ \to N^\circ$ gives the surjection $M^\circ \to \Hom(N_\uf\cap N^\circ,\Z)$ with kernel $N_\uf^\perp$, and we have the identification $\Hom(N_\uf\cap N^\circ,\Z)\simeq M^\circ/N_\uf^\perp$ for $p_2$.  We choose a map $p:N\to M^\circ$ compatible with both $p_1$ and $p_2$ by restriction or quotient.
\end{defn}
Let $[a]_+:=\max(0,a)$ and $[a]_-:=\min(0,a)$. 
\begin{defn}Given a seed $\s$ and $k\in I_\uf$, we define the \emph{mutation} $\mu_k:N\to N$ by
\begin{equation}
\mu_k(\bfe_i)=
\begin{cases}
\bfe_i+[\ve_{ki}]_+\bfe_k &i\neq k\\
-\bfe_k & i=k
\end{cases}.
\end{equation}
\end{defn}

It is obvious that $\{\mu_k(\bfe_i)\}_{i\in I}$ is still a basis for $N$, $\{\mu_k(\bfe_i)\}_{i\in I_\uf}$ is still a basis for $N_{\uf}$, and $\{d_i\mu_k(\bfe_i)\}$ is still a basis for $N^\circ$. Then for the new seed $\mu_k(\s):=\{\mu_k(\bfe_i)\}_{i\in I}$, one checks that
\begin{equation}
\mu_k(\bff_i):=d_i^{-1}\mu_k(\bfe_i)^*=\begin{cases}
\bff_i &i\neq k\\
-\bff_k+\sum_j [-\ve_{jk}]_+\bff_j & i=k
\end{cases},
\end{equation}
and the exchange matrix changes as
\begin{equation}
\label{ExchMatrixMutation}
\mu_k(\ve_{ij}):=d_i\{\mu_k(\bfe_i),\mu_k(\bfe_j)\}=
\begin{cases}
-\ve_{ij} & k\in\{i,j\}\\
\ve_{ij} & \ve_{ik}\ve_{kj}\leq 0, k\not\in\{i,j\}\\
\ve_{ij}+|\ve_{ik}|\ve_{kj} & \ve_{ik}\ve_{kj}> 0, k\not\in\{i,j\}
\end{cases}.
\end{equation}
Note that this matrix mutation coincides with the mutation of the exchange matrix defined in \cite{FZ02}. It can also be written as
\begin{equation}
\label{ExchMatrixMutation2}
\mu_k(\ve_{ij})=
\begin{cases}
-\ve_{ij} & k\in\{i,j\}\\
\ve_{ij}+\ve_{ik}[\ve_{kj}]_+ + [-\ve_{ik}]_+\ve_{kj} & k\not\in\{i,j\}
\end{cases}.
\end{equation}

Let $\fT$ be the infinite oriented rooted tree with $|I_\uf|$ outgoing edges from each vertex, labeled by the elements of $I_\uf$. We will label the root of $\fT$ by $t_0$, and we attach a seed $\s_0$ to $t_0$, called the \emph{initial seed}. Each node $t\in \fT$ corresponds to a sequence of mutations and a seed $\s_t$ obtained from the initial seed. We will use subscript $t\in\fT$ to specify which seed various terms are obtained from, and we might omit the subscript for $t=t_0$ if we do not need to emphasize the seed. For example, if $t_0 \xrightarrow{k} t_1$, then $\ve_{ij;t_1}:=\mu_k(\ve_{ij})$.\\
 
Let $\bbk$ be a field of characteristic zero. For any lattice $N$ with $M=\Hom(N,\Z)$, the algebraic torus $T_N$ over $\bbk$ is the scheme $T_N:=\Spec \bbk[M]$ where $\bbk[M]$ is isomorphic to the group ring consisting of $\bbk$-linear combinations of $z^m$ for $m\in M$ where $z$ is a formal variable.

For a given seed $\s$, we consider the algebraic tori 
\Eq{\cX_\s&:=T_{M,\s}=\Spec \bbk[N]\\
\cA_\s&:=T_{N^\circ,\s}=\Spec \bbk[M^\circ],}
where the subscript $\s$ indicates the choice of seed. 
\begin{defn}For $k\in I_{\uf}$ we define the birational maps
\begin{align*}
\mu_k&: \cX_\s\rat \cX_{\mu_k(\s)}\\
\mu_k&: \cA_\s \rat \cA_{\mu_k(\s)}
\end{align*}
using pullback of functions
\begin{align}
\label{z^nMutation}
\mu_k^*(z^n)&=z^n(1+z^{\bfe_k})^{-\{d_k\bfe_k,n\}}\qquad n\in N    \\
\label{z^mMutation}
\mu_k^*(z^m)&=z^m(1+z^{\bfv_k})^{-\langle d_k\bfe_k,m\rangle} \qquad m\in M^\circ
\end{align}
where $\bfv_k=p(\bfe_k)=\{\cdot,\bfe_k\}\in M^\circ$ and $\langle\cdot,\cdot\rangle$ is the pairing between $N^\circ$ and $M^\circ$. Note that these maps implicitly depend on the choice of $\s=(\bfe_i)_{i\in I}$.
\end{defn}

For $X_i:=z^{\bfe_i}$, $A_i:=z^{\bff_i}$, $X_i':=z^{\mu_k(\bfe_i)}$, and $A_i':=z^{\mu_k(\bff_i)}$, the mutation formulas are equivalent to
\begin{equation}
\label{XMutation}
\mu_k^*(X_i')=\begin{cases}
X_iX_k^{[\ve_{ki}]_+}(1+X_k)^{-\ve_{ki}} & i\neq k\\
X_k^{-1} & i=k
\end{cases}
\end{equation}
and
\begin{equation}
\label{AMutation}
\mu_k^*(A_i')=\begin{cases}
A_i & i\neq k\\
A_i^{-1}\left(\prod_j A_j^{[\ve_{jk}]_+}+\prod_j A_j^{[-\ve_{jk}]_+}\right) & i=k
\end{cases}.
\end{equation}
As explained in the beginning, these formulas differ from the mutation formulas for cluster variables and universal rational functions defined in \cite{FZ02} and \cite{FZ07} by a transpose of $\ve_{jk}$.

\begin{defn}
The \textit{Fock-Goncharov $\cA$, $\cX$ cluster varieties} is obtained by gluing the $\cA_\s$ and $\cX_\s$ over all seeds $\s$ mutation equivalent to an initial seed $\s_0$ using the mutation birational maps.
\end{defn}

The map $p:N\to M^\circ$ induces a map $p:\cA\to \cX$ given by $p^*(z^n)=z^{p(n)}$, by directly checking
\begin{align*}
p^*(\mu_k^*(z^n))&=z^{p(n)}(1+z^{p(\bfe_k)})^{-\{d_k\bfe_k,n\}},\\
\mu_k^*(p^*(z^n))&=z^{p(n)}(1+z^{\bfv_k})^{-\langle d_k\bfe_k,p(n)\rangle},
\end{align*}
so $p^*(\mu_k^*(z^n))=\mu_k^*(p^*(z^n))$ by the definition of $p$. In terms of the $A$ and $X$ variables, the map $p$ is given by
\begin{equation}
p^*(X_i)=\A\prod_{j\in I_\uf} A_j^{\ve_{ji}},
\end{equation}
where $\A$ is a product of $A_k$ over frozen $k\in I\setminus I_\uf$, which depends on the choice of $p$.\\

The Laurent phenomenon is a well-known result on the regularity of functions on $\cA$ and $\cX$. For a proof, see \cite{FZ02}, \cite[Section 3]{GHK13}.

\begin{thm}[Laurent phenomenon]
\label{LaurentPhenomenon}
For any seed $\s$,
\begin{enumerate}
\item Let $z^m\in \bbk[M^\circ]$ such that $\langle m,d_i\bfe_i\rangle\geq 0$ for each $i\in I_\uf$. Then $z^m\in \cO(\cA)$. 
\item Let $z^n\in \bbk[N]$ such that $\{\bfe_i,n\}\geq 0$ for each $i\in I_\uf$. Then $z^n\in \cO(\cX)$.
\end{enumerate}
In particular, any cluster $\cA$-variable $A_i:=z^{\bff_i}$ is a regular function on $\cA$.
\end{thm}

Such functions are called \textit{global monomials} in \cite{GHKK}. In other literature, they are called \textit{cluster monomials} for $\cA$ and \textit{standard monomials} for $\cX$ respectively.\\\\

Finally, we recall the tropicalization $\cX(\Z^\trop)$ as follows. Let $\Z^\trop$ be the semifield $(\Z,+,\min)$, and $Q_{sf}(M)$ be the set of subtraction-free rational functions in the field of fractions of $\bbk[N]$, which is a semifield. We define $T_M(\Z^\trop):=\Hom_{sf}(Q_{sf}(M),\Z^\trop)$. Since the monomials $\{z^n\}_{n\in N}$ generate $Q_{sf}(M)$, we get a bijection
\[
\Hom_{sf}(Q_{sf}(M),\Z^\trop)\to \Hom_{group}(N,\Z)=M.
\]
By the mutation formula, each mutation $\mu_k^*$ on $\bbk[N]$ is invertible and involves only subtraction-free expressions, so it induces an isomorphism on $Q_{sf}(M)$. This allows us to glue the $T_M(\Z^\trop)$ to give $\cX(\Z^\trop)$. Each choice of vertex in $\fT$ gives an identification $\cX(\Z^\trop)=M$.

We fix using an initial seed $\s_0$ for this identification, which fixes a choice of $\bfe_i^*\in \cX(\Z^\trop)$. For each $t\in \fT$, let
\Eq{\label{tropx}
x_{ij;t}:=\bfe_i^*(\mu_t^*(X_{j;t})) \in\Z^\trop,
}
where $\mu_t$ is the composition of pullbacks determined by $t$. Then for $t\xrightarrow{k} t'$ in $\fT$, applying $\bfe_i^*$ to the mutation formula gives
\begin{equation}
\label{TropMutation}
x_{ij,t'}=\begin{cases}
-x_{ik;t} & \text{if $j=k$}\\
x_{ij;t}+[\ve_{kj;t}]_+x_{ik;t}+\ve_{kj;t}[-x_{ik;t}]_+ & \text{if $j\neq k$}
\end{cases}.
\end{equation}

We also give another way of tropicalizing. Instead of fixing the initial seed and $\bfe_i^*\in \cX(\Z^\trop)$, we let
\[
y_{ij;t}:=\bfe_{i;t}^*((\mu_t^{-1})^*(X_j))\in\Z^\trop.
\]
Here $\bfe_{i;t}^*$ changes with $t$, but we fix the cluster monomials $X_j$ in the initial seed.

\begin{remark}
\label{NotationIssue}
On the level of cluster $\cA$- and $\cX$-variables, the definitions used here and the definitions in \cite{FG2} and \cite{GHKK} are related by taking a transpose of the exchange matrix $(\ve_{ij})$ introduced in detail in \cite[Section 1]{FG2}. However, for the study of $G$-local systems, the two notations give exactly the same cluster algebras because of the different choice of the Cartan matrix, which is also off by a transpose.
\end{remark}

\subsection{$c$ and $g$-vectors}
\label{sec:prelim:cvector}

In this section, we introduce the $c$ and $g$-vectors, which were originally introduced in \cite{FZ07} to describe the mutation of cluster algebras with principal coefficients. Here we follow the reinterpretation in \cite{GHKK}, which allowed them to prove the sign-coherence conjectures in the skew-symmetrizable case.\\

We first need to define the principal fixed data, which constructs the cluster variety with principal coefficients.

\begin{defn}
\label{PrincipalFixedData}
Given a fixed data $\Gamma$, the \textit{principal fixed data} $\Gamma_\prin$ is the fixed data with:
\begin{itemize}
\item $\wt{N}=N\oplus M^\circ$ with the $\Q$-valued skew-symmetric bilinear form
\Eq{
\{(n_1,m_1),(n_2,m_2)\}=\{n_1,n_2\}-\langle n_1,m_2\rangle+\langle n_2,m_1\rangle.
}
Here the pairing $M^\circ \x N^\circ\to \Z$ is be extended to $M^\circ \x N\to \Q$ by $\langle \bfe_i,m\rangle=d_i^{-1}\langle d_i\bfe_i,m\rangle$.
\item $\wt{N}_\uf=N_\uf\oplus 0\subset \wt{N}$.
\item $\wt{N}^\circ=N^\circ\oplus M$.
\item $\wt{M}=M\oplus N^\circ$, $\wt{M}^\circ=M^\circ\oplus N$.
\item The index set $\wt{I}$ is the disjoint union of two copies of $I$, with $d_i$ the same as in $\Gamma$. The unfrozen indices is the original $I_\uf$ in the first copy of $I$.
\end{itemize}
\end{defn}

\begin{defn}
For any choice of initial seed $\s=(\bfe_1,\cdots,\bfe_n)$ for $\Gamma$, the initial seed of $\wt{\s}$ for $\Gamma_\prin$ is chosen to be
\Eq{
\wt{\s}:=((\bfe_1,0),\cdots,(\bfe_n,0),(0,\bff_1),\cdots,(0,\bff_n))
}
and gives a mutation class $[\wt{\s}]$. For each $t\in \fT$, we will denote by $\wt{\s}_t$ the mutation of $\wt{\s}$ along the sequence of unfrozen indices defined by $t$.
\end{defn}

We will denote the corresponding $\cA$ and $\cX$ cluster varieties for $\Gamma_\prin$ by $\cA_\prin$ and $\cX_\prin$ respectively.

Note that for $\s_t=(\bfe_{1;t},\cdots,\bfe_{n;t})$, it is not true that
\[
\wt{\s}_t=((\bfe_{1;t},0),\cdots,(\bfe_{n;t},0),(0,\bff_{1;t}),\cdots,(0,\bff_{n;t})).
\]
However, by definition of the mutation rule the first half is still correct. Thus we still have a lift $A_{i;t}=z^{\bff_{i;t}}$ from $\cA_{\s_t}\subset \cA$ to $\cA_{\prin,\s_t}\subset \cA_\prin$, where we identify $\bff_{i;t}$ with $(\bff_{i;t},0)\in \wt{M}^\circ$.

We will denote the indices in the second copy of $I$ by $i'$, corresponding to each $i$ in the first copy of $I$. In terms of the quiver $Q$, this corresponds to adding a frozen vertex $i'$ for each $i$ in $Q$ and an arrow $i'\to i$ since $\{(\bfe_i,0),(0,\bff_i)\}=-1$ by definition.

\begin{defn}
For each $t\in \fT$, the \textit{$C$-matrix} $C_t$ for $\s_t$ is the $n\x n$-matrix with entries
\Eq{
c_{ij}=\wt{\ve}_{i'j;t}:=d_i\{\bfe_{i'},\bfe_j\},\quad i,j\in I
}
The column vectors $\bfc_{j;t}$ of $C_t$ are called the \textit{$c$-vectors}.
\end{defn}

For $g$-vectors, consider the surjection $\wt{M}^\circ\to M^\circ$ given by $(m,n)\mapsto m-p(n)$. This gives torus action of $T_{N^\circ}$ on $T_{\wt{N}^\circ}$. Explicity, let $m-p(n)=\sum_{i\in I} c_i\bff_i\in M^\circ$. Then a point $(a_i)\in T_{N^\circ}$ acts by
\[
(a_i)\cdot z^{(m,n)}=\prod_{i\in I} a_i^{c_i} z^{(m,n)}.
\]
This action is compatible with mutation \ref{z^mMutation} since we have $\{\cdot,(\bfe_k,0)\}=(p(\bfe_k),\bfe_k)$, so $z^{\{\cdot,(\bfe_k,0)\}}$ has weight 0. Thus we get a torus action of $T_{N^\circ}$ on $\cA_\prin$. Clearly any cluster monomial $z^{(m,n)}$ is a $T_{N^\circ}$-eigenfunction under this action with weight $m-p(n)$.

\begin{defn}
The \textit{$g$-vector} $\bfg_{j;t}$ is the $T_{N^\circ}$-weight of the lift of $A_{j;t}$ to $\cA_\prin$. It is identified with a vector in $\Z^n$ using the basis $(\bff_1,\cdots,\bff_n)$.
\end{defn}

The following theorem on the sign coherence of $c$ and $g$-vectors is conjectured in \cite[Section 6]{FZ07} and proved in \cite[Section 5]{GHKK}.
\begin{thm}
\label{SignCoherence}
For each $j\in I$, the entries of $\bfc_{j;t}$ are either all non-negative or all non-positive. For each $i\in I$, the $i$-th coordinate of $\bfg_{j;t}$ for any $j\in I$ are either all non-negative or all non-positive.
\end{thm}

For later use, it is convenient to have an explicit mutation formula for the $c$ and $g$-vectors, which is originally given in \cite{FZ07}.

\begin{lem}
The $c$-vectors are given by the initial condition
\Eq{
c_{ij;t_0}=\delta_{ij}
}
and the mutation formula
\begin{equation}
\label{CMutation}
c_{ij;t'}=\begin{cases}
-c_{ik;t} & \text{if $j=k$}\\
c_{ij}+c_{ik}[\ve_{kj}]_+ + [-c_{ik}]_+\ve_{kj} & \text{if $j\neq k$}
\end{cases}
\end{equation}
for $t \xrightarrow{k} t'$ in $\fT$.
\end{lem}

\begin{proof}
Immediate from the mutation formula in \ref{ExchMatrixMutation} and the initial conditions in \ref{PrincipalFixedData}.
\end{proof}

\begin{lem}
The $g$-vectors are given by the initial condition
\Eq{
g_{ij;t_0}=\delta_{ij}
}
and the mutation formula
\begin{equation}
\label{GMutation}
g_{ij;t'}=\begin{cases}
g_{ij;t} & \text{if $j\neq k$}\\
-g_{ik;t}+ \sum_{l\in I} g_{il;t}[\ve_{lk;t}]_+ - \sum_{l\in I} \ve_{il;t_0}[c_{lk;t}]_+ & \text{if $j=k$}
\end{cases}
\end{equation}
for $t \xrightarrow{k} t'$ in $\fT$.
\end{lem}

\begin{proof}
By the mutation formula in \ref{AMutation},
\[
\mu_k^*(A_{j;t'})=A_{j;t}^{-1}\left(\prod_{l\in I} A_{l;t}^{[\ve_{lk;t}]_+}\prod_{l\in I}A_{l';t}^{[\ve_{l'k;t}]_+}+\prod_{l\in I} A_{l;t}^{[-\ve_{lk;t}]_+}\prod_{l\in I}A_{l';t}^{[-\ve_{l'k;t}]_+}\right).
\]
Note that $A_{l';t}=A_{l';t_0}$ is fixed under mutation because $l'$ is frozen.

The two monomials in the mutation formula above must have the same $T_{N^\circ}$-weight because $\mu_k^*(A_{j;t'})$ and each $A_{l;t}$, $A_{l';t}$ are $T_{N^\circ}$-eigenvectors. Thus we can calculate the weight using the first term. By definition $c_{lk;t}=\ve_{l'k;t}$, and the weight of $A_{l';t_0}=z^{(0,\bfe_l)}$ is $-p(\bfe_l)$. Hence the $T_{N^\circ}$-weight of $A_{j;t'}$ is given by
\[
\bfg_{j;t'}=-\bfg_{j;t}+\sum_{l\in I}[\ve_{lk;t}]_+\bfg_{l;t}-\sum_{l\in I} [c_{lk;t}]_+p(\bfe_{l;t_0}).
\]
On the root $t_0\in\fT$, the $\bff_i$ coefficient of $p(\bfe_l)$ is
\[
d_i\langle \bfe_i,p(\bfe_l)\rangle=d_i\{\bfe_i,\bfe_l\}=\ve_{il},
\]
and the result follows.
\end{proof}

\begin{remark}
By calculating using the second term instead, the mutation formula \ref{GMutation} can also be written as
\begin{equation}
\label{GMutation2}
g_{ij;t'}=\begin{cases}
g_{ij;t} & \text{if $j\neq k$}\\
-g_{ik;t}+ \sum_{l\in I} g_{il;t}[-\ve_{lk;t}]_+ - \sum_{l\in I} \ve_{il;t_0}[-c_{lk;t}]_+ & \text{if $j=k$}
\end{cases}.
\end{equation}
\end{remark}

Finally, we note that the mutation of $c$-vectors in \ref{CMutation} agrees with the tropicalized mutation \ref{TropMutation}. Thus we get the following lemma
\begin{lem}
\label{C=X}
For any $t\in \fT$, $x_{ij;t}=c_{ij;t}$.
\end{lem}

Meanwhile, using \cite[Eq. (1.13)]{NZ12} and Theorem \ref{SignCoherence}, we get the following lemma
\begin{lem}
\label{G=Y}
Let $g_{ij}^\dagger$ be the $g$-vectors associated to the transposed seed defined in Remark \ref{NotationIssue}. Then for any $t\in \fT$, 
\Eq{y_{ij;t}=g^\dagger_{ji;t}.}
\end{lem}

\subsection{Quantization}
\label{sec:prelim:quantum}

To define quantum cluster algebras, we quantize the definitions in Section \ref{sec:prelim:defn} as follows:

\begin{defn}
\label{Quantization}
Let $q$ be a formal parameter. For any finite rank lattice $N$ with $\frac{1}{d}\Z$-valued skew-symmetric bilinear form $\{\cdot,\cdot\}$, the \textit{quantum torus algebra} $\bbk_q[N]$ is defined to be the noncommutative algebra over $\bbk_q=\bbk[q^{\pm 1/d}]$ generated by $z^n$ for $n\in N$, subject to the relations
\[
q^{\{n_1,n_2\}}z^{n_1}z^{n_2} = z^{n_1+n_2}.
\]

For $\cO_q(\cX_\s):=\bbk_q[N]$ and $\cO_q(\cA_\s):=\bbk_q[M^\circ]$, the quantization of the mutations \ref{z^nMutation} and \ref{z^mMutation} are given by the pullback maps, and the image are elements on the division ring (well defined by the Ore condition):
\begin{align}
\label{z^nMutationQuantum}
\mu_k^*(z^n)&:=\Ad_{\Psi_{q_k}(z^{\bfe_k})}z^n,\qquad n\in N\\
\label{z^mMutationQuantum}
\mu_k^*(z^m)&:=\Ad_{\Psi_{q_k}(z^{\bfv_k})}z^m\qquad m\in M^\circ,
\end{align}
where $q_k:=q^{1/d_k}$, and $\Psi_{q_k}(z^{\bfe_k})$ is the \emph{quantum dilogarithm} defined by
\begin{equation}
\label{QuantumDilog}
\Psi_q(x):=\prod_{r=1}^\infty \frac{1}{1+q^{2r-1}x},
\end{equation}
while $\Ad_\Psi$ is the conjugation by $\Psi$.

In parallel to the definition of $\cX$ and $\cA$, we define the \textit{quantum upper cluster algebras} $\cO_q(\cX)$ and $\cO_q(\cA)$ to be
\begin{align}
\cO_q(\cX)&:=\bigcap_{t\in \fT} \mu_t^*(\cO_q(\cX_{\s_t})),\\
\cO_q(\cA)&:=\bigcap_{t\in \fT} \mu_t^*(\cO_q(\cA_{\s_t})),
\end{align}
which contain rational functions that are regular in any seed.
\end{defn}

For $X_i=z^{\bfe_i}$ and $X_i'=z^{\bfe_i'}$, a direct computation shows that \ref{z^nMutationQuantum} can be rewritten as
\begin{equation}
\label{XMutationQuantum}
\mu_k^*(X_i')=\begin{cases}
X_k^{-1} & i=k\\
X_i\prod_{r=1}^{|\ve_{ki}|}(1+q_k^{2r-1}X_k) & i\neq k,\ve_{ik}\leq 0\\
X_i\prod_{r=1}^{|\ve_{ki}|}(1+q_k^{2r-1}X_k^{-1})^{-1} & i\neq k,\ve_{ik}\geq 0
\end{cases},
\end{equation}
while for $A_i'=z^{\bff_i'}$, \ref{z^mMutationQuantum} can be rewritten as
\begin{equation}
\label{AMutationQuantum}
\mu_k^*(A_i')=\begin{cases}
A_i & i\neq k\\
z^{\bff_k'}+z^{\bff_k'+\bfv_k} & i=k
\end{cases},
\end{equation}
Note that taking $q=1$ we recover the classical mutation formulas \ref{XMutation} and \ref{AMutation}.

For later use, we will write
\[
X_{i_1,\cdots,i_k}:=z^{\bfe_{i_1}+\cdots+\bfe_{i_k}}
\]
for $i_1,\cdots,i_k\in I$. Note that $X_{i_1,\cdots,i_k}\neq X_{i_1}\cdots X_{i_k}$.\\

Many of the good properties in classical cluster algebras can be extended to quantum cluster algebras. For example, we still have the injection $p^*:\cO_q(\cX)\to \cO_q(\cA)$. The quantum Laurent phenomenon is proved in \cite{BZ05} for $A$ variables, and is extended to $X$ variables using $p^*$ in \cite{GS}.

\begin{thm}[Quantum Laurent phenomenon]
\label{QuantumLaurentPhenomenon}
For any seed $\s$,
\begin{enumerate}
\item Let $z^v\in \bbk_q[M^\circ]$ such that $\langle v,e_i\rangle\geq 0$ for each $i\in I_\uf$. Then $z^v\in \cO_q(\cA)$. 
\item Let $z^v\in \bbk_q[N]$ such that $\{e_i,v\}\geq 0$ for each $i\in I_\uf$. Then $z^v\in \cO_q(\cX)$.
\end{enumerate}
In particular, any quantum cluster $\cA$-variable $A_i=z^{\bff_i}$ belongs to $\cO_q(\cA_q)$.
\end{thm}

The tropicalization introduced in Section \ref{sec:prelim:defn} can be extended to $\cO_q(\cX)$ by restricting to the semi-field homomorphisms that send $q^{\pm 1/D}$ to 0. Then it is checked easily that the tropical mutation for quantum $X$ variables agree exactly with the classical ones.

\subsection{Amalgamation}
\label{sec:prelim:amal}

Another important construction of cluster algebras is amalgamation, which corresponds to the gluing of two quivers along frozen variables. This construction introduced in \cite{FG3} appears naturally in the study of cluster algebras arising from higher Teichm\"uller theory, since we glue quivers associated to triangles along their boundaries to get quivers for a general surface.

\begin{defn}
Let $\{(\Gamma_a,\s_a)\}_{a\in A}$ be a finite collection of seed datum. For any index set $I$ with surjection
\[
\rho: \bigsqcup_{a\in A} I_a\to I
\]
such that:
\begin{itemize}
\item The restriction $\rho_a:=\rho|_{I_a}$ is injective for all $a\in A$,
\item If $\rho(j)=\rho(k)$ for $j\in I_a$ and $k\in I_{a'}$, then $j,k$ are frozen and $d_j=d_k$.
\end{itemize}
The \textit{amalgamation} is a new seed data $(\Gamma,\s)$ where
\begin{itemize}
\item $N$ is the sublattice of $\bigoplus_{a\in A}N_a$ spanned by $\s=\{\bfe_i\}_{i\in I}$, where
\[
\bfe_i=\sum_{\rho(j)=i} \bfe_j,
\]
\item The bilinear form $\{\cdot,\cdot\}$ is induced from the direct sum $\bigoplus_{a\in A}N_a$.
\item For each $i\in I$, $d_i=d_j$ for any $j\in \rho^{-1}(i)$.
\end{itemize}
\end{defn}

It is easy to see that the cluster algebra $\cO(\cA_{\Gamma})$ is a subalgebra of the tensor product $\bigotimes_{a\in A}\cO(\cA_{\Gamma_a})$, and we have natural injections $\cO(\cA_{\Gamma_a})\to \cO(\cA_\Gamma)$ given by $A_j\mapsto A_{\rho(j)}$. The same construction holds on the cluster $\cX$-side as well as their quantizations.

\section{Quantum groups and positive representations}
\label{sec:repn}

\subsection{Root systems}
\label{sec:repn:root}

Let $\fg$ be a finite dimensional semisimple Lie algebra over $\C$ of rank $r$. Let $I=\{1,\cdots,r\}$ be the indexing set for the Dynkin diagram of $\fg$. Let $\Phi$ be the root system of $\fg$ with inner product $(\cdot,\cdot)$, $\Pi_+=\{\alpha_i\}_{i\in I}$ be a set of simple positive roots, and $\Phi_+$ be the set of positive roots. Let
\Eq{
a_{ij}:=\frac{2(\alpha_i,\alpha_j)}{(\alpha_i,\alpha_i)},\qquad i,j\in I
}
and $A:=(a_{ij})$ be the Cartan matrix of $\fg$. Each simple positive root comes with an integer multiplier $d_i$ with
\Eq{\label{di}
d_i^{-1}:=\frac{1}{2}(\alpha_i,\alpha_i)
}
such that $d_i^{-1}a_{ij}=d_j^{-1}a_{ji}$.\footnote{This definition of $d_i$ is the inverse of the one defined in \cite{Ip18} to match the multipliers in the fixed data. The parameter $q_i$ agrees so the definitions of the quantum cluster algebras and quantum groups are the same as in \cite{Ip18}.} We have $d_i\in \{1,2,3\}$.

We also have the simple positive coroots $\{\alpha_i^\vee\}_{i\in I}$ such that $\langle \alpha_i,\alpha_j^\vee\rangle = a_{ij}$, and the fundamental weights $\{\Lambda_i\}_{i\in I}$ such that $\langle \Lambda_j,\alpha_i^\vee\rangle = \delta_{ij}$.

The Weyl group $W$ is the group of isometries of $\Phi$ generated by the simple reflections $s_i:=s_{\alpha_i}$. It has a unique longest element $w_0$, which induces an involution $*:I\to I$ such that
\begin{equation}
\alpha_{i^*}=-w_0(\alpha_i).
\end{equation}

\subsection{Quantum groups $\cU_q(\fg)$ and $\fD_q(\fg)$}
\label{sec:repn:Uq}

Let $q_i=q^{1/d_i}$, same as in Section \ref{sec:prelim:quantum}, and $d=lcm\{d_i\}_{i\in I}$. For $n,m\in \Z_{\geq0}$, define
\begin{align*}
[n]_q&:=\frac{q^n-q^{-n}}{q-q^{-1}}, & [n]_q!&:=\prod_{s=1}^n[s]_q, & \begin{bmatrix} n\\ m\\ \end{bmatrix}_q:=\frac{[n]_q!}{[m]_q![n-m]_q!}.
\end{align*}

\begin{defn}
We define $\fD_q(\fg)$ to be the $\C(q^{1/d})$-algebra generated by $\{\bfE_i,\bfF_i,\bfK_i^{\pm1},{\bfK_i'}^{\pm1}\}_{i\in I}$ with the following relations:
\begin{align}
\bfK_i\bfK_j&=\bfK_j\bfK_i, & \bfK_i'\bfK_j'&=\bfK_j'\bfK_i',\\
\bfK_i\bfK_j'&=\bfK_j'\bfK_i, & [\bfE_i,\bfF_j]&=\delta_{ij}(q_i-q_i^{-1})(\bfK_i'-\bfK_i),\\
\bfK_i\bfE_j&=q^{a_{ij}}\bfE_j\bfK_i, & \bfK_i\bfF_j&=q^{-a_{ij}}\bfF_j\bfK_i,\\
\bfK_i'\bfE_j&=q^{-a_{ij}}\bfE_j\bfK_i', & \bfK_i'\bfF_j&=q^{a_{ij}}\bfF_j\bfK_i',
\end{align}
and also the quantum Serre relations
\begin{align}
\sum_{s=0}^{1-a_{ij}} (-1)^s \begin{bmatrix} 1-a_{ij}\\ s \end{bmatrix}_q \bfE_i^k\bfE_j\bfE_i^{1-a_{ij}-k}&=0,\\
\sum_{s=0}^{1-a_{ij}} (-1)^s \begin{bmatrix} 1-a_{ij}\\ s \end{bmatrix}_q \bfF_i^k\bfF_j\bfF_i^{1-a_{ij}-k}&=0.
\end{align}
\end{defn}

\begin{defn}
$\fD_q(\fg)$ has a Hopf algebra structure, with coproduct
\begin{align}
\Delta(\bfE_i)&=1\ox \bfE_i + \bfE_i\ox \bfK_i, & \Delta(\bfK_i)&=\bfK_i\ox \bfK_i,\\
\Delta(\bfF_i)&=\bfF_i\ox 1 + \bfK_i'\ox \bfF_i, & \Delta(\bfK_i')&=\bfK_i'\ox \bfK_i',
\end{align}
counit
\begin{align}
\ve(\bfE_i)&=\ve(\bfF_i)=0, & \ve(\bfK_i)&=\ve(\bfK_i')=1,
\end{align}
and antipode
\begin{align}
S(\bfE_i)&=-\bfK_i^{-1}\bfE_i, & S(\bfK_i)&=\bfK_i^{-1},\\
S(\bfF_i)&=-\bfF_i\bfK_i, & S(\bfK_i')&={\bfK_i'}^{-1}.
\end{align}
\end{defn}

\begin{defn}
The quantum group $\cU_q(\fg)$ is given by the quotient
\begin{equation}
\cU_q(\fg)=\fD_q(\fg)/\langle \bfK_i\bfK_i'=1\rangle_{i\in I}.
\end{equation}
\end{defn}
Note that $\fD_q(\fg)$ can be identified with the \emph{Drinfeld's double} of the Borel part $\cU_q(\fb)$ of the quantum group.

\subsection{Cluster realizations of $\fD_q(\fg)$}
\label{sec:repn:cluster}

\subsubsection{Punctured disk realization}
\label{sec:repn:cluster:punctured}

In \cite{Ip18}, the first author constructed a quiver $\bfD(\bfi_0)$ for any reduced word $\bfi_0$ of $w_0$, and showed that there is an embedding $\fD_q(\fg)$ into the quantum cluster algebra $\cO_q$ for this quiver. The image of $\bfF_i$ are telescoping sums along certain $F_i$-paths, and the image of the $\bfK_i'$ generators are monomials along the $F_i$-paths. When $\bfi_0$ is appropriately chosen, for example if $i_N=i$ where $N=l(w_0)$, so are the images of $\bfE_i$ and $\bfK_i'$.

\begin{thm}[{\cite[Theorem 4.14]{Ip18}}]
\label{IpEmbedding}
Given a reduced word $\bfi_0$ of the longest element of the Weyl group, there is a quiver $\bfD(\bfi_0)$ of rank $2l(w_0)+2r$ giving a quantum torus algebra $\cO_q(\cX_{\bfD(\bfi_0)})$, with an embedding
\begin{equation}
\label{IpDqEmbedding}
\iota: \fD_q(\fg)\hookrightarrow \cO_q(\cX_{q,\bfD(\bfi_0)}),
\end{equation}
which descends to an embedding
\begin{equation}
\label{IpUqEmbedding}
\iota: \cU_q(\fg)\hookrightarrow \cO_q(\cX_{q,\bfD(\bfi_0)})/\langle \iota(\bfK_i)\iota(\bfK_i')=1\rangle_{i\in I}.
\end{equation}
\end{thm}

These embeddings a priori depend on the choice of reduced word $\bfi_0$. In \cite[Section 7]{Ip18}, it is also shown that for different choices of $\bfi_0$, the corresponding quivers $\bfD(\bfi_0)$ are related by mutations for each braid relation. Based on this observation, the first author suspected that there are embeddings
\begin{align}
\fD_q(\fg)&\hookrightarrow \cO_q(\cX), & \cU_q(\fg)&\hookrightarrow \cO_q(\cX)/\cI,
\end{align}
which would be independent of $\bfi_0$.

This conjecture is proved in \cite{GS}. Goncharov and Shen constructed the moduli space of $G$-local systems $\sP_{G,\bbS}$, where $G$ is a split semisimple adjoint group over $\Q$ and $\bbS$ is a decorated oriented surface satisfying some mild conditions. A $G$-local system in $\sP_{G,\bbS}$ is a principal $G$-bundle $\cL$ on $\bbS$ with flat connection with the following extra data:
\begin{enumerate}
\item At boundary special point $m$, a flat section of the associated decorated flag bundle $\cL_\cA=\cL\x_G\cA$ on a small neighborhood of $m$.
\item At each puncture $p$, a flat section of the associated flag bundle $\cL_\cB=\cL\x_G\cB$ on a small neighborhood of $p$.
\end{enumerate}

It is shown in \cite[Section 10]{GS} that $\sP_{G,\bbS}$ carries a cluster Poisson structure, i.e. corresponds to a cluster $\cX$-variety. There is an atlas for each choice of ideal triangulation for $\bbS$, along with a choice of vertex and also a choice of reduced word of $w_0$ for each triangle of the ideal triangulation. Each of these atlases are related by a sequence of mutations. We then quantize this cluster algebra and obtain the quantum upper cluster algebra $\cO_q(\sP_{G,\bbS})$.

There are functions $\cW_{s,i}, \cK_{s,i}\in \cO_q(\sP_{G,\bbS})$ for each boundary special point $s$ and $i\in I$, constructed geometrically from the $G$-local system definition of $\sP_{G,\bbS}$, then lifted to the quantum cluster algebra. The functions $\cW_{s,i}$ are called \emph{potential functions}, which are conjectured to be the Landau-Ginzburg potentials in mirror symmetry.

When $\bbS=\odot$ is the once punctured disk with two boundary special points $s,t$, the quiver for $\sP_{G,\odot}$ is mutation equivalent to the quiver constructed in \cite{Ip18}. There is an outer monodromy map
\[
\mu_{out}:\sP_{G,\odot}\to H,
\]
such that the condition $\mu_{out}=1$ corresponds to $\cK_{s,i}\cK_{t,i^*}=1$ for all $i\in I$. By \cite[Lemma 15.5]{GS}, the functions $\cW_{s,i}$ and $\cK_{s,i}$ become monomials after a sequence of mutations, which satisfy the assumptions in Theorem \ref{QuantumLaurentPhenomenon}. Therefore, $\cW_{s,i}, \cK_{s,i}\in \cO_q(\sP_{G,\odot})$, and we have the following maps into the usual quantum upper cluster algebra $\cO_q(\sP_{G,\odot})$ and also $\cO_q(\sP_{G,\odot})_{\mu_{out}=1}$, improving Theorem \ref{IpEmbedding}.

\begin{thm}[{\cite[Theorem 3.21]{GS}}]
\label{GSEmbedding}
There is a canonical embedding
\begin{equation}
\label{GSDqEmbedding}
\kappa: \fD_q(\fg)\to \cO_q(\sP_{G,\odot})
\end{equation}
given by $\bfE_i\mapsto \cW_{s,i}$, $\bfF_i\mapsto \cW_{t,i^*}$, $\bfK_i\mapsto \cK_{s,i}$, $\bfK_i'\mapsto \cK_{t,i^*}$, which descends to an embedding
\begin{equation}
\label{GSUqEmbedding}
\kappa: \cU_q(\fg)\to \cO_q(\sP_{G,\odot})_{\mu_{out}=1}.
\end{equation}
\end{thm}

The image is in fact invariant under a Weyl group action at the puncture defined in \cite[Section 13]{GS}, but we do not need it and we will not describe it in detail here.

In \cite{She22}, it is shown that when $\fg$ is simply laced, the embedding $\kappa:\cU_q(\fg)\to \cO_q(\sP_{G,\odot})$ descends to an isomorphism
\[
\cU_q(\fg)\simeq \cO_q(\sL_{G,\odot}),
\]
where $\sL_{G,\odot}$ is obtained from $\sP_{G,\odot}$ by first forgetting the framing at the punctures giving the space $\textnormal{Loc}_{G,\odot}$, then taking the fiber $\mu_{out}^{-1}(1)$ over the identity $1\in H$. This is also conjectured to be true for general $\fg$. Note that the map forgetting the framing
\[
\sP_{G,\odot}\to \textnormal{Loc}_{G,\odot}
\]
is a Galois $W$-cover, so the map $\kappa$ descends to $\cO_q(\textnormal{Loc}_{G,\odot})$ exactly by the fact that the image of $\cU_q(\fg)$ is invariant under the Weyl group action. 

\subsubsection{Parabolic realization}
\label{sec:repn:cluster:parabolic}

Based on Theorem \ref{IpEmbedding}, the first author constructed a cluster realization of $\fD_q(\fg)$ for each parabolic subgroup $P_J\subset G$ for any subset $J\subset I$ of Dynkin nodes \cite{Ip20}. The simple reflections $(s_j)_{j\in J}$ generate the parabolic subgroup $W_J\subset W$ of the Weyl group. Let $w_J$ be the longest element of $W_J$. Then we can decompose $w_0$ as
\[
w_0=w_J\overline{w}
\]
for a unique $\overline{w}\in W$, and any reduced word $\overline{\bfi}$ for $\overline{w}$ gives a quiver $\bfD(\overline{\bfi})$.\\

Let $\cO_q(\overline{\cX})$ be the quantum cluster algebra associated to $\bfD(\overline{\bfi})$.
\begin{thm}[{\cite[Theorem 5.2]{Ip20}}]
There is a homomorphism
\begin{equation}
\overline{\kappa}: \fD_q(\fg) \to \cO_q(\overline{\cX}).
\end{equation}
\end{thm}

\subsubsection{Zero Casimirs realization}
\label{sec:repn:cluster:zero}

In \cite{SS17}, Schrader and Shapiro described a mutation sequence for $\bfD(\bfi)$, where $\fg=\mathfrak{sl}_{n+1}$ and
\begin{equation}
\bfi=(1,2,1,3,2,1,\cdots ,n,n-1,\cdots,1),
\end{equation}
such that the new quiver contains a symmetric part with $n+1$ variables. Based on this mutation sequence, the first author and Man performed folding of these variables and obtained a new quiver and new quantum cluster algebra $\cO_q(\cX^0)$ \cite{IM22}. In view of the embedding in \ref{IpDqEmbedding}, the folding of cluster variables are related to the Casimir operators $\bfC_k$ in $\fD_q(\fg)$, which gives the following result:
\begin{thm}[{\cite[Theorem 4.5]{IM22}}]
There is an embedding
\begin{equation}
\label{ZeroCasimirEmbedding}
\kappa_0:\fD_q(\mathfrak{sl}_{n+1})/\langle \bfC_k=0\rangle \to \cO_q(\cX^0).
\end{equation}
\end{thm}

Using an inclusion of a parabolic $A_k$ part for general Lie types $\fg$, the embedding \ref{ZeroCasimirEmbedding} can be generalized as follows:
\begin{thm}[{\cite[Theorem 5.1]{IM22}}]
For any parabolic subgroup $P_J\subset G$ of type $A_k$, there is a homomorphism
\begin{equation}
\label{ZeroCasimirHomomorphism}
\kappa_0:\fD_q(\fg) \to \cO_q(\cX^0).
\end{equation}
where $\cO_q(\cX^0)$ is obtained by folding the symmetric part of the parabolic $A_k$ part of the quiver $\bfD(\bfi_0)$.
\end{thm}

\section{Polynomial conditions}
\label{sec:polynomial}

\subsection{Main results}

To prepare for the main theorem, we first discuss the sufficient conditions for monomials or polynomials in the cluster $\cX$-variables for one seed to be universally monomial or universally polynomial. Consider a fixed data $\Gamma$ with a fixed initial seed $\s_0=(\bfe_i)_{i\in I}$ as in Definition \ref{fixeddata}--\ref{defseed}.

\begin{lem}
\label{PrelimMonomial}
Let $z^v\in\bbk_q[N]$ be such that $\{\bfe_i,v\}=0$ for each $i\in I_\uf$. Then $z^v$ is always a monomial in the quantum cluster variables for any seed in the mutation class of $\s_0$.
\end{lem}

\begin{proof}
Since $z^v$ and $z^{-v}$ both satisfy the assumptions in Theorem \ref{QuantumLaurentPhenomenon}, they are Laurent polynomials in the cluster variables for any seed in the mutation class, hence they must be monomials.
\end{proof}

\begin{prop}
\label{PrelimPolynomial}
For any seed $\s_t$, $t\in \fT$, let
\[
v=\sum_{i\in I-I_\uf} a_i \bfe_{i;t},
\] 
for some $a_i\in \Z_{\geq0}$, and suppose $\{\bfe_{j;t},v\}\geq 0$ for each $j\in I_\uf$. Then $z^v\in \bbk_q[N]$ remains a polynomial in the quantum cluster variables for any seed in the mutation class.
\end{prop}

\begin{proof}
It suffices to check at the initial seed. By Theorem \ref{QuantumLaurentPhenomenon}, $\mu_t^*(z^v)$ is a Laurent polynomial.

Meanwhile, consider the tropicalization \eqref{tropx}. For each $l\in I$, using Lemma \ref{C=X} we get
\[
\bfe_l^*(\mu_t^*(z^v))=\sum_{i\in I-I_\uf} a_ix_{li;t}=\sum_{i\in I-I_\uf} a_ic_{li;t}.
\]
Since the frozen vertices $i\in I-I_\uf$ are never mutated, $c_{ii;t}=1$ by definition. By sign coherence of $c$-vectors, $c_{li;t}\geq 0$ for all $l\in I$. This shows that $\bfe_l^*(\mu_t^*(z^v))\geq 0$ for all $l\in I$, so the power of each $X_l$ in $\mu_t^*(z^v)$ is non-negative and $\mu_t^*(z^v)$ is a polynomial.
\end{proof}

The condition in Proposition \ref{PrelimPolynomial} is sufficient but not necessary. For example, for the quantum cluster algebra associated to the Markov quiver (without frozen vertices), the monomial $X_1X_2X_3$ is invariant under mutation, hence universally polynomial.\\

It is interesting to ask exactly which standard monomials are universally polynomial, or more generally, what exactly is the subring of universal polynomials in $\cO(\cX)$ or $\cO_q(\cX)$. A necessary condition is given as follows.

Fix an initial seed $\s_0$ with the corresponding cluster $\cX$-variables $X_1,\cdots,X_n$. Let $f\in \cO_q(\cX)$ be a polynomial in $X_i$, such that each monomial appearing in $f$ determines a vector $\vec{v}$ in $\Z^n_{\geq0}$. Let $\wt{G}_t$ be the matrix for $g_{ji;t}^\dagger$ as described in Lemma \ref{G=Y}. 
\begin{prop}
\label{PrelimPolyIff} $f$ is universally polynomial if and only if $\wt{G}_t\vec{v}\in \Z^n_{\geq0}$ for all $t\in \fT$ and all $\vec{v}$.
\end{prop}

\begin{proof}
The forward direction is proved in the same as Proposition \ref{PrelimPolynomial}. Conversely, suppose $\wt{G}_t\vec{v}\not\in \Z^n_{\geq0}$ for some $t\in\fT$ and $\vec{v}$. By \cite[Proposition 4.2]{NZ12} and Theorem \ref{SignCoherence}, $\wt{G}_t$ is invertible, so $\wt{G}_t\vec{v}$ is distinct for different $\vec{v}$. Recall from Lemma \ref{G=Y} and the definition of $y_{ij;t}$ that $\wt{G}_t\vec{v}$ corresponds to the lowest degree term with respect to $\s_t$. This shows that its negative powers do not cancel and $f$ has a term with negative powers in the seed $\s_t$.
\end{proof}

While Proposition \ref{PrelimPolyIff} characterizes all universal polynomials, we do not know a good way to compute the exact condition, as there are in general infinitely many $\wt{G}_t$ when $t$ ranges through $\fT$. However, using the language of \cite{GHKK}, $g$-vectors are generators of the chambers in the cluster complex $\Delta^+(\Z)$, which is the collection of cluster chambers, each of which is a cone in the tropicalization $\cA^\vee(\Z^\trop)$. Thus we can reinterpret the condition as follows.

\begin{cor}
Let $f\in \cO_q(\cX)$. Then $f$ is universally polynomial if and only if for each monomial $z^n$ appearing in $f$ in some seed $\s$, $\langle n,m\rangle \geq0$ for each $m\in \Delta^+(\Z)$.
\end{cor}

\begin{eg}
If $\cA$ has a large cluster complex, $\Delta^+(\Z)$ does not lie in a half space, so there are no non-trivial universal polynomials in $\cO_q(\cX)$. This happens if for example $\cA$ has an acyclic seed, or has a seed with a maximal green sequence \cite[Corollary 8.30]{GHKK}.

This shows that in many cases, including the moduli space of local systems $\cP_{G,\bbS}$ in \cite[Theorem 2.13]{GS}, frozen variables are necessary for $\cO_q(\cX)$ to have non-trivial universal polynomials.
\end{eg}

\begin{eg}
Consider the standard quiver described in Section \ref{sec:repn:cluster} for the embedding of $\fD_q(\sl_2)$ in Figure \ref{fig-A0}, where $1,3$ are frozen. 

\begin{figure}[htb!]
\centering
\begin{tikzpicture}[every node/.style={inner sep=0, minimum size=0.5cm, thick, fill=white, draw}, x=2cm, y=1cm]
\node (1) at (0,1) {$1$};
\node (2) at (1,0) [circle]{$2$};
\node (3) at (2,1) {$3$};
\node (4) at (1,2)  [circle]{$4$};
\drawpath{1,2,3}{black}
\drawpath{3,4,1}{black}
\end{tikzpicture}
\caption{The quiver for $\cX_q^{\mathrm{std}}$.}\label{fig-A0}
\end{figure}

It is easy to compute that
\Eq{
\cO(\cX)=\bbk[X_2^\pm, X_4^\pm, (X_1X_3)^\pm, X_3(1+X_4),X_1(1+X_2)].
}
This quiver is skew-symmetric, so the transposed seed is just itself. The matrices $\wt{G}_t$ for the 4 possible seeds are
\[\begin{pmatrix}
1 & 0 & 0 & 0\\
0 & 1 & 0 & 0 \\
0 & 0 & 1 & 0 \\
0 & 0 & 0 & 1 \\
\end{pmatrix},
\begin{pmatrix}
1 & 0 & 0 & 0\\
1 & -1 & 0 & 0 \\
0 & 0 & 1 & 0 \\
0 & 0 & 0 & 1 \\
\end{pmatrix},
\begin{pmatrix}
1 & 0 & 0 & 0\\
0 & 1 & 0 & 0 \\
0 & 0 & 1 & 0 \\
0 & 0 & 1 & -1 \\
\end{pmatrix},
\begin{pmatrix}
1 & 0 & 0 & 0\\
1 & -1 & 0 & 0 \\
0 & 0 & 1 & 0 \\
0 & 0 & 1 & -1 \\
\end{pmatrix}\]
respectively, so the condition in Proposition \ref{PrelimPolyIff} is just $d_1\geq d_2$ and $d_3\geq d_4$ for each monomial $X_1^{d_1}X_2^{d_2}X_3^{d_3}X_4^{d_4}$ with all $d_i\geq 0$. Then the subring of universal polynomials is given explicitly by
\Eq{
\bbk[X_1(1+X_2), X_3(1+X_4), X_1X_3,X_1X_2X_3, X_1X_3X_4,X_1X_2X_3X_4],\label{sl2uni}
}
which can also be written as (using $'$ to denote the variables in mutated seeds)
\[
\bbk[X_1X_3,\mu_2^*(X_1'),\mu_2^*(X_1'X_3'),\mu_4^*(X_3''),\mu_4^*(X_1''X_3''),\mu_2^*\mu_4^*(X_1'''X_3''')].
\]
Thus in this case the monomials in Proposition \ref{PrelimPolynomial} generate the universal polynomials.
\end{eg}

\begin{remark}
In Proposition \ref{PrelimPolyIff}, the criterion does not determine whether $f\in \bbk_q[X_1,\cdots,X_n]$ belongs to $\cO_q(\cX)$, so we need to assume $f\in \cO_q(\cX)$. It raises the question on the exact size of $\cO(\cX)$ or $\cO_q(\cX)$, for example whether it is generated by standard monomials. This question on the $\cX$-side parallels the question whether a cluster algebra coincides with its upper cluster algebra, which happens on the $\cA$-side \cite{BFZ05}.
\end{remark}

The following lemma gives a class of functions which are universally polynomial.

\begin{lem}
\label{AnCase}
Let $\cX$ be a cluster variety defined by a quiver $Q$. Suppose $Q$ contains a type $A_n$ subquiver $Q'$ indexed by $1,\dots,n$, where the vertices $1$ and $n$ are frozen and there is an arrow $(i\to i+1)$ for each $i=1,...,n-1$. Suppose for any unfrozen vertex $v\in Q-Q'$,
\begin{itemize}
\item there are no arrows $(1\to v)$,
\item for any $v'\in Q'$, there is at most one arrow $(v\to v')$ or $(v'\to v)$,
\item whenever there is an arrow $v\to v'$, there is another $v''\in Q'$ with $(v''>v')$ with an arrow $(v''\to v)$, and $v$ is non-adjacent to any vertices in $Q'$ between $v'$ and $v''$.
\end{itemize}
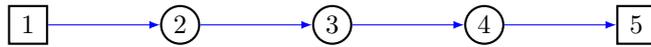
\begin{figure}[htb!]
\centering
\begin{tikzpicture}[every node/.style={inner sep=0, minimum size=0.5cm, thick, fill=white, draw}, x=2cm, y=1cm]
\node (1) at (0,0) {$1$};
\node (2) at (1,0) [circle]{$2$};
\node (3) at (2,0) [circle]{$3$};
\node (4) at (3,0) [circle]{$4$};
\node (5) at (4,0) {$5$};
\drawpath{1,2,3,4,5}{blue}
\end{tikzpicture}
\caption{The type $A_n$ quiver for $n=5$.}\label{fig-A1}
\end{figure}
Then the polynomials
\begin{align}
X_1+X_{1,2}+\cdots+X_{1,2,\cdots,n-1}&:=z^{\bfe_1}+z^{\bfe_1+\bfe_2}+\cdots+z^{\bfe_1+\bfe_2+\cdots+\bfe_{n-1}}\label{X123}\\
X_{1,2,\cdots,n}&:=z^{\bfe_1+\bfe_2+\cdots+\bfe_n}\label{K123}
\end{align}
are universally polynomial.
\end{lem}

\begin{proof}
Starting from Figure \ref{fig-A1}, we first mutate at $2$ to get the quiver in Figure \ref{fig-A2}.

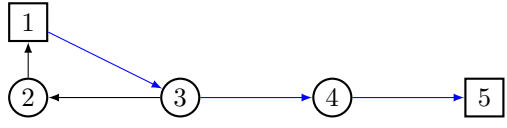
\begin{figure}[htb!]
\centering
\begin{tikzpicture}[every node/.style={inner sep=0, minimum size=0.5cm, thick, fill=white, draw}, x=2cm, y=1cm]
\node (1) at (1,1) {$1$};
\node (2) at (1,0) [circle]{$2$};
\node (3) at (2,0) [circle]{$3$};
\node (4) at (3,0) [circle]{$4$};
\node (5) at (4,0) {$5$};
\drawpath{1,3,4,5}{blue}
\drawpath{3,2,1}{black}
\end{tikzpicture}
\caption{The quiver in Figure \ref{fig-A1} after mutating at $2$.}\label{fig-A2}
\end{figure}

By direct computation,
\begin{align*}
\mu_2^*(X_1'+X_{1,3}'+X_{1,3,4}'+\cdots+X_{1,3,4,\cdots,n-1}')&=X_1+X_{1,2}+\cdots+X_{1,2,\cdots,n-1},\\
\mu_2^*(X_{1,3,4,\cdots,n}')&=X_{1,2,\cdots,n}.
\end{align*}

Then we mutate at $3,\cdots,n-1$ to get the quiver in Figure \ref{fig-A3}.

\begin{figure}[htb!]
\centering
\begin{tikzpicture}[every node/.style={inner sep=0, minimum size=0.5cm, thick, fill=white, draw}, x=2cm, y=1cm]
\node (1) at (3,1) {$1$};
\node (2) at (1,0) [circle]{$2$};
\node (3) at (2,0) [circle]{$3$};
\node (4) at (3,0) [circle]{$4$};
\node (5) at (4,0) {$5$};
\drawpath{1,5}{blue}
\drawpath{5,4,1}{black}
\drawpath{2,3,4}{black}
\end{tikzpicture}
\caption{The quiver in Figure \ref{fig-A1} after mutating at $2,\cdots,n-1$.}\label{fig-A3}
\end{figure}
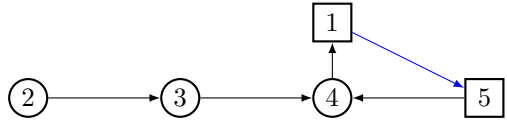

Inductively, we get the result
\begin{align*}
\mu^*(X_1'')&=X_1+X_{1,2}+\cdots+X_{1,2,\cdots,n-1},\\
\mu^*(X_{1,n}'')&=X_{1,2,\cdots,n}.
\end{align*}

In order to apply Proposition \ref{PrelimPolynomial}, we need to descibe the arrows between $1,n$ and vertices in $Q-Q'$. For any $v\in Q-Q'$, initially if there are arrows between $v$ and 1, it must be a single arrow $(v\to 1)$, then by assumption there is another $v'\in Q'$ with an arrow $(v'\to v)$. Then the arrow $(v\to 1)$ cancels after mutating at $(2,3,\dots,v')$. Similarly, all arrows between $v$ and $n$ cancel after mutation. Thus the only arrows involving $1$ and $n$ at the last step of the mutations are $(n\to n-1)$ and $(n-1\to n)$ only as seen in Figure \ref{fig-A3}, and the assumption in Proposition \ref{PrelimPolynomial} is satisfied.
\end{proof}

Finally, we introduce and prove the main theorem of this paper.

\begin{thm}
\label{MainThm}
For the homomorphisms
\begin{align*}
\kappa&: \fD_q(\fg) \to \cO_q(\cX)\\
\overline{\kappa}&: \fD_q(\fg) \to \cO_q(\overline{\cX})\\
\kappa_0&:\fD_q(\fg) \to \cO_q(\cX^0)
\end{align*}
described in Section \ref{sec:repn:cluster}, the image of the Chevalley generators $\bfE_i,\bfF_i,\bfK_i,\bfK_i'$ are all polynomials in the quantum cluster variables with coefficients in $\N[q^{\pm1}]$ for any cluster in the mutation class.
\end{thm}

\begin{proof}
The image of $\bfF_i$ and $\bfK_i'$ in $\cO_q(\cX)$ are given by
\begin{align*}
\bfF_i&\mapsto X_{{i\choose -n_i}}+X_{{i\choose -n_i},{i\choose -n_i+1}}+\cdots+X_{{i\choose -n_i},\cdots,{i\choose n_i-1}},\\
\bfK_i'&\mapsto X_{{i\choose -n_i},{i\choose -n_i+1},\cdots,{i\choose n_i-1}}.
\end{align*}
Here ${i\choose k}$ are the indices used in \cite{GS}, where each $i$ corresponds to a positive simple root $\alpha_i$, and $n_i$ is the number of times $\alpha_i$ appears in the reduced expression of $w_0$ used when choosing the seed for $\cX$. By construction, the quiver $\bfD(\bfi)$ is amalgamated from triangular quivers based on the Dynkin diagram, thus the assumptions in Lemma \ref{AnCase} is satisfied, and $\bfF_i,\bfK_i'$ are mapped to polynomials in the cluster variables for any cluster in the mutation class.\\

For $\bfE_i$ and $\bfK_i$, we can choose a reduced word $\bfi$ for $w_0$ such that $i_N=i$. Then it is known \cite{Ip18} that the image of $\bfE_i$ and $\bfK_i$ are mutation equivalent to the telescopic sum and monomial respectively as in \eqref{X123}--\eqref{K123}, along the $E_i$-path, which is just an $A_3$ quiver. Lemma \ref{AnCase} again shows that $\bfE_i,\bfK_i$ are mapped to polynomials in the cluster variables for any cluster in the mutation class.\\

For the parabolic representations, the argument for $\bfF_i$ and $\bfK_i'$ is the same since the $F_i$ paths remain type $A_n$. For $\bfE_i$ and $\bfK_i$, we might no longer be able to choose freely the reduced word $i_N=i$ due to restriction by the decomposition $w_0=w_J\overline{w}$. We instead use the mutations involved in the proof of \cite[Lemma 5.18]{Ip20}, which shows that the image of $\bfE_i,\bfF_i,\bfK_i,\bfK_i'$ are regular. It is shown that there is a sequence of quantum cluster mutations that send the images into a monomial (a single variable in the case of $\bfE_i$ and $\bfF_i$) having the same adjacency to the quiver of $\bfD(\bfi)$ in the standard embedding $\kappa$. Therefore, the images also satisfy the conditions in Proposition \ref{PrelimPolynomial}, and the same conclusion follows for $\cO_q(\overline{\cX})$.\\

For the zero Casimir representations, the result follows immediately from a general property of cluster folding Corollary \ref{FoldingPolynomial}, which we will discuss in Section \ref{sec:polynomial:folding}.\\

As for the coefficients, the result follows from the main theorem of \cite{DM21}, where Davison and Mandel showed that for skew-symmetric quantum cluster algebras, the quantum theta functions have coefficients in $\N[q^{\pm1}]$ in any seed. Since standard monomials are theta functions, this immediately implies our result when $G$ is simply laced. For general $G$, we can unfold the cluster algebra to obtain a skew-symmetric one. Standard monomials in the folded cluster algebra come from standard monomials in the unfolded one, so they also have coefficients in $\N[q^{\pm1}]$ in any seed.
\end{proof}

\begin{eg}
For any simple $\fg$, the generalized Casimir elements \cite{KS98} are given by
\Eq{
\textbf{C}_k:=(1\ox \tr|^q_{V_k})(RR_{21})
}
where
\begin{itemize}
\item $V_k$ is the $k$-th fundamental representation of $\cU_q(\fg)$ for $k=1,\cdots,n$,
\item $\tr|^q_{V_k}$ is the quantum trace over $V_k$, and
\item $R$ is the universal $R$-matrix.
\end{itemize}
The generalized Casimir elements lie in $\widehat{\cU}_q(\fg):=\cU_q(\fg)[\bfK_i^{\pm\frac{1}{h}}]$, where $h$ is the Coxeter number of $\fg$, and they generate the center of $\widehat{\cU}_q(\fg)$. We can modify this construction to $\fD_q(\fg)$ as well, and it can be shown that in $\widehat{\fD}_q(\fg)$, 
\Eq{
\textbf{C}_k=\cK_k\cC_k
}
for each $k$, where $\cK_k$ is the Cartan part consists of product of rational powers of $\bfK_i,\bfK_i'$, while $\cC_k$ is a polynomial in $\bfE_i,\bfF_i,\bfK_i,\bfK_i'$. Then Theorem \ref{MainThm} shows that the polynomial part is universally polynomial in $\cO_q(\cX)$.\\

For example, in $\mathfrak{sl}_2$, the unique generalized Casimir element is given by
\Eq{
\textbf{C}=\bfF\bfE-q\bfK-q^{-1}\bfK'.
}
There is no Cartan part here, and $\textbf{C}$ is indeed a polynomial in $\bfE_i,\bfF_i,\bfK_i,\bfK_i'$.

In the standard embedding $\fD_q(\mathfrak{sl}_2)\hookrightarrow\cO_q(\cX)$, the image of $\textbf{C}$ is $X_{1,2,3,4}+X_{1,3}$ which is indeed universally polynomial as shown in \eqref{sl2uni}.
\end{eg}

\subsection{Cluster folding}
\label{sec:polynomial:folding}

The positive representations with zero Casimirs involves cluster folding, originally defined in \cite[Version 3 Section 6]{FG2} and \cite[Section 3.6]{FG3} to apply folding of root systems to the moduli spaces of $G$-local systems. To apply Theorem \ref{MainThm}, we will define a slightly more restrictive version requiring more symmetry to obtain better properties.

\begin{defn}
\label{FoldingDefn}
For any fixed data $\Gamma$ with seed $\s=(\bfe_i)_{i\in I}$, suppose there is a finite group $\Sigma$ acting on $I$ preserving $I_\uf\subset I$ such that
\Eq{
\{\bfe_{\sigma(i)},\cdot\}=\{\bfe_i,\cdot\}
}
for all $i\in I$ and $\sigma\in \Sigma$. Clearly $\Sigma$ induces an automorphism of $N$.
Then we define the \textit{folding} of $\Gamma$ to be a fixed data $\Gamma^\Sigma$ consisting of:
\begin{itemize}
\item The index set $\overline{I}:=I/\Sigma$, $\overline{I}_\uf:=I_\uf/\Sigma$, $d_{\overline{i}}:=d_i/|O(i)|$, where $O(i)$ is the orbit of $i\in I$ under $\Sigma$.
\item $N^\Sigma$ is the invariant sublattice generated by $\overline{\s}=(\bfe_{\overline{i}})_{\overline{i}\in \overline{I}}$ where
\Eq{
\bfe_{\overline{i}}:=\sum_{j\in O(i)}\bfe_j,
}
with $N_\uf^\Sigma$ generated by $(\bfe_{\overline{i}})_{\overline{i}\in \overline{I}_\uf}$, and $(N^\Sigma)^\circ$ generated by $(d_{\overline{i}}\bfe_{\overline{i}})_{\overline{i}\in \overline{I}}$.
\item The skew-symmetric bilinear form on $N^\Sigma$ is inherited from that on $N$, so
\Eq{
\{\bfe_{\overline{i}},\bfe_{\overline{j}}\}:=|O(i)||O(j)|\{\bfe_i,\bfe_j\}.
}
\item Lattices $M^\Sigma:=\Hom(N^\Sigma,\Z)$, $(M^\Sigma)^\circ:=\Hom((N^\Sigma)^\circ,\Z)$.
\end{itemize}
\end{defn}

From the definition, it is obvious that the requirements $\{(N^\Sigma)^\circ,N_\uf^\Sigma\}\subset \Z$, $\{N_\uf^\Sigma\cap (N^\Sigma)^\circ,N^\Sigma\}\subset \Z$ are always satisfied. The numbers $d_{\overline{i}}$ might not be integers, but this does not affect any of our calculations and we can always fix this by rescaling if necessary.\\

In the remainder of this section, we will always assume that $\Gamma$ has a folding $\overline{\Gamma}$ defined using the initial seed $\s_0$.\\

Define analogously $\bff_{\overline{i}}:=d_{\overline{i}}^{-1} \bfe_{\overline{i}}^*$, $\ve_{\overline{i}\overline{j}}:=d_{\overline{i}}\{\bfe_{\overline{i}},\bfe_{\overline{j}}\}$ and the cluster varieties $\cX^\Sigma:=\cX_{\Gamma^\Sigma}$, $\cA^\Sigma:=\cA_{\Gamma^\Sigma}$ and similarly $\cO_q(\cX^\Sigma)$, $\cO_q(\cA^\Sigma)$ respectively.

Let $\overline{\fT}$ be the infinite oriented rooted tree for $\Gamma^\Sigma$. We define an injection $\iota:\overline{\fT}\to \fT$ as follows:
\begin{enumerate}
\item The root $\overline{t}_0$ of $\overline{\fT}$ is mapped to the root $t_0$ of $\fT$.
\item Inductively, suppose $\iota(t)$ is already defined and $t\xrightarrow{\overline{k}} t'$ in $\overline{\fT}$. We choose an arbitrary order $k_1,\cdots,k_{|O(k)|}$ of the orbit $O(k)$ of $k$ and move from $\iota(t)$ along the edges labeled by $k_1,\cdots,k_{|O(k)|}$ in such order giving the vertex $\iota(t')$.
\end{enumerate}

The following lemma shows that the folding $\Gamma^\Sigma$ is independent of the seed $\s_t$ for $t\in \iota(\overline{\fT})$.

\begin{lem}
For any $\overline{t}\in \overline{\fT}$, the seed $\s_{\iota(\overline{t})}$ also defines a folding. Then $N^\Sigma_{\overline{t}}=N^\Sigma$, and the following diagram commutes
\begin{center}
\begin{tikzcd}
N^\Sigma \ar[r]\ar[d,"\mu_{\overline{t}}"] & N \ar[d,"\mu_{\iota(\overline{t})}"]\\
N^\Sigma \ar[r] & N
\end{tikzcd}
\end{center}
where the horizontal maps are inclusions.
\end{lem}

\begin{proof}
We reduce to the case $\overline{t}_0\xrightarrow{\overline{k}}\overline{t}$ by induction.

By the assumption in Definition \ref{FoldingDefn} and the mutation formulas in Section \ref{sec:prelim:defn}, if $i$ and $j$ have the same orbit, the mutations $\mu_i$ and $\mu_j$ are identical up to swapping $i$ and $j$. In this case, $\ve_{ij}=0$, and it is also well known that the mutations $\mu_i$ and $\mu_j$ commute. Thus we still have
\[
\{\mu_i\mu_j(\bfe_i),\cdot\}=\{\mu_i\mu_j(\bfe_j),\cdot\}.
\]
This shows that the assumption in Definition \ref{FoldingDefn} is still satisfied and inductively $\s_{\iota(\overline{t})}$ also defines a folding.

The condition $N^\Sigma_{\overline{t}}=N^\Sigma$ is immediate since mutation sends $N^\Sigma$ into itself. Let $\mu_{O(k)}$ be the composition of the commuting mutations $\mu_l$ over all $l\in O(k)$. Then
\[
\mu_{O(k)}(\bfe_i)=
\begin{cases}
\bfe_i+[\ve_{ki}]_+\bfe_k &i\not\in O(k)\\
-\bfe_i & i\in O(k)
\end{cases},
\]
so we have
\[
\mu_{O(k)}\circ(\bfe_{\overline{i}})=
\begin{cases}
\bfe_{\overline{i}}+|O_i|[\ve_{ki}]_+\bfe_{\overline{k}} &i\not\in O(k)\\
-\bfe_{\overline{i}} & i\in O(k)
\end{cases}.
\]
Since
\[
|O(i)|\ve_{ki}=|O(i)|d_k\{\bfe_k,\bfe_i\}=\ve_{\overline{k}\overline{i}},
\]
we get $\mu_{\overline{t}}=\mu_{\iota(\overline{t})}$ on $N^G$.
\end{proof}

The inclusion $N^\Sigma\to N$ also induces an injection $\iota:\bbk_q[N^\Sigma]\to \bbk_q[N]$. Unlike for $N$ above, the mutations for $\cX$ does not commute with $\iota$, so we cannot glue them to an injection $\cO_q(\cX^\Sigma)\to \cO_q(\cX)$ globally.

However, we note that the image of $\iota$ lies the the invariant subring $\bbk_q[N]^\Sigma$. It is well known that $\bbk_q[N]^\Sigma$ is generated over $\bbk_q$ by $\cE_{i,a}$ for $i\in I$, $|a|\leq |O(i)|$, where $\cE_{i,a}$ is the elementary symmetric polynomial of degree $|a|$ in the variables $X_j$ for $a>0$ and $X_j^{-1}$ for $a<0$, over all $j\in O(i)$. We have $\iota(X_{\overline{i}})=\cE_{i,|O(i)|}$.

Then by extending the map $\iota$ and the mutation maps to the completion of $\bbk_q[N^G]$ and $\bbk_q[N]^\Sigma$ with respect to the ideal generated by the cluster variables $X_i$, we can still salvage the following fact:

\begin{prop}
The injection $\iota:\Frac(\widehat{\bbk_q[N^\Sigma]})\to \Frac(\widehat{\bbk_q[N]^\Sigma})$ has a $\bbk_q$-linear inverse $\pi$ such that the following diagram commutes
\begin{center}
\begin{tikzcd}
\Frac(\widehat{\bbk_q[N]^\Sigma}) \ar[r,"\pi"]\ar[d,"\mu_{\iota(\overline{t})}"] & \Frac(\widehat{\bbk_q[N^\Sigma]}) \ar[d,"\mu_{\overline{}}"]\\
\Frac(\widehat{\bbk_q[N]^\Sigma}) \ar[r,"\pi"] & \Frac(\widehat{\bbk_q[N^\Sigma]})
\end{tikzcd}
\end{center}
\end{prop}

\begin{proof}
For any $j \in O(i)$ and $k\not\in O(i)$ with $\ve_{ki}\leq 0$, we have
\[
X_kX_j=q^{-2\{\bfe_k,\bfe_j\}}X_jX_k=q_k^{2|\ve_{ki}|}X_jX_k.
\]
This shows that for $i_1,\cdots,i_a\in O(i)$,
\begin{align*}
\mu_k^*(X_{i_1}'\cdots X_{i_a}')&= X_{i_1}\prod_{r=1}^{|\ve_{ki}|}(1+q_k^{2r-1}X_k)\cdots X_{i_a}\prod_{r=1}^{|\ve_{ik}|}(1+q_k^{2r-1}X_k)\\
&= X_{i_1}\cdots X_{i_a}\prod_{r=1}^{|\ve_{ki}|}\prod_{s=0}^{a-1} (1+q_k^{2r+2s|\ve_{ik}|-1}X_k)\\
&= X_{i_1}\cdots X_{i_a}\prod_{r=1}^{a|\ve_{ki}|} (1+q_k^{2r-1}X_k).
\end{align*}
We also have similar formulas when $\ve_{ik}\geq0$.

Then the mutation $\mu_{O(k)}^*$ on the elementary symmetric polynomials are given by
\begin{equation}
\label{SymMutation}
\mu_{O(k)}^*(\cE_{i,a}')=\begin{cases}
\cE_{i,-a}^{-1} & O(i)= O(k)\\
\cE_{i,a}\prod_{r=1}^{a|\ve_{ki}|}(\sum_{b=0}^{|O(k)|} q_k^{(2r-1)b}\cE_{k,b}) & O(i)\neq O(k),\ve_{ki}\leq 0\\
\cE_{i,a}\prod_{r=1}^{a|\ve_{ki}|}(\sum_{b=0}^{|O(k)|} q_k^{(2r-1)b}\cE_{k,-b})^{-1} & O(i)\neq O(k),\ve_{ki}\geq 0
\end{cases}.
\end{equation}

Meanwhile,
\begin{equation}
\label{SymMutation}
\mu_{\overline{k}}^*(X_{\overline{i}}')=\begin{cases}
X_{\overline{k}}^{-1} & \overline{i}=\overline{k}\\
X_{\overline{i}}\prod_{r=1}^{|O(i)||\ve_{ki}|}(1+q_{\overline{k}}^{2r-1}X_{\overline{k}}) & \overline{i}\neq \overline{k},\ve_{\overline{i}\overline{k}}\leq 0\\
X_{\overline{i}}\prod_{r=1}^{|O(i)||\ve_{ki}|}(1+q_{\overline{k}}^{2r-1}X_{\overline{k}}^{-1})^{-1} & \overline{i}\neq \overline{k},\ve_{\overline{i}\overline{k}}\geq 0
\end{cases}.
\end{equation}

Note that $q_{\overline{k}}=q_k^{|O(i)|}$. A direct comparison shows that $\iota$ identifies the mutation of $X_{\overline{i}}$ with the terms of the mutation of $\cE_{i,|O(i)|}$ in $\Frac(\widehat{\bbk_q[N]^G})$ involving only $\cE_{j,|O(j)|}$ for $j\in I$. This shows that the $\bbk_q$-linear map $\pi$ sending $\cE_{i,b|O(i)|}$ to $X_{\overline{i}}^b$ and all other $\cE_{i,a}$ to 0 is the inverse of $\iota$ and commutes with mutation. 
\end{proof}

We immediately get the following corollary which completes the proof of Theorem \ref{MainThm}.

\begin{cor}
\label{FoldingPolynomial}
Suppose $f\in \bbk_q[N]$ is a polynomial in the quantum cluster variables for any seed in the mutation class of $\s$. Then so is $\pi(f)$ for any seed in the mutation class of $\overline{\s}$.
\end{cor}

\subsection{Extension to general $G$-local systems}
\label{sec:polynomial:tensor}

Our result can be extended to embeddings of $\fD_q(\fg)$ into $\cO_q(\sP_{G,\bbS})$ for general $G$-local systems. The main example is $\odot_n$, the disk with $n$-punctures and two boundary special points. The relation between quantum groups and multiple punctured disk is first investigated in \cite{Kas01} for $n=2$, \cite{SS19} for general $n$ and $\fg=\mathfrak{sl}_{r+1}$, and in \cite{GS} for general $\fg$.

Geometrically, $\odot_n$ is obtained by gluing $n$ copies of $\odot$ along boundary intervals, thus by amalgamation, the cluster algebra $\cO_q(\sP_{G,\odot_n})$ is a subalgebra of the tensor product $\cO_q(\sP_{G,\odot})^{\ox n}$. For example, one checks that for $n=2$,
\begin{align}
\cW_{s,i,\odot_2}&=\cW_{s,i,\odot}\ox 1+\cK_{s,i,\odot}\ox \cW_{s,i,\odot}\\
\cK_{s,i,\odot_2}&=\cK_{s,i,\odot}\ox \cK_{s,i,\odot}
\end{align}
and the formula for general $n$ can be derived inductively. Combining with the embedding $\kappa$, we obtain the formula for the coproduct in $\cU_q(\fg)$.

\begin{cor}
\label{PolynomialGeneralization}
$\cW_{s,i,\odot_n}$ and $\cK_{s,i,\odot_n}$ remain polynomials in the quantum cluster variables with coefficients in $\N[q^{\pm1}]$ for any cluster in the mutation class.
\end{cor}

\begin{proof}
For each $i$, the cluster $\cX$-variables appearing in $\cW_{s,i,\odot_n}$ and $\cK_{s,i,\odot_n}$ are still in the configuration of a type $A_n$ quiver. Therefore, the polynomiality result follows from Lemma \ref{AnCase}. The statement about the coefficients again follows from the main theorem of \cite{DM21}.
\end{proof}

\begin{remark}
Note that the proof for the parabolic case in \ref{MainThm} also relies on mutations into polynomials only involving frozen variables, and folding can be applied to each $\cO_q(\sP_{G,\odot})$ separately. Therefore, the same result still holds if for some copies of $\odot$, the quiver for $\cO_q(\sP_{G,\odot})$ is replaced by the parabolic or zero Casimir variant.
\end{remark}

For general marked surfaces, we no longer have the symmetry in $\odot$ or $\odot_n$, so for any boundary special point $s$, the embedding is restricted to only
\[
\cU_q(\fb)\to \cO_q(\sP_{G,\bbS}),
\]
where the choice of Borel subalgebra $\fb$ is determined by the Borel subgroup from the framing at $s$ \cite[Theorem 3.21]{GS}. Lemma \ref{AnCase} again shows that the image of $\bfE_i$ and $\bfK_i$ are also universally polynomial in this case.

As mentioned in \cite[Remark 15.3]{GS}, the above result can be generalized to arbitrary Kac-Moody $\fb$, where the map is no longer injective and the RHS is defined directly using a quiver without using $G$-local systems. Our polynomiality result extends naturally as well. Furthermore, our result also applies to the rare example \cite[Proposition 4.10]{Ip20} which maps the whole quantum affine algebra $\cU_q(\widehat{\mathfrak{sl}}_{n+1})$ to certain quantum cluster algebra via the evaluation module.

\section{Regular positive representations}
\label{sec:discussion}
After proving Theorem \ref{MainThm}, we would like to discuss another conjecture proposed in \cite[Section 6]{IM22}. The first author and Man introduced the notion of \textit{regular positive representations}, which are positive representations $\cP$ of $\cU_q(\fg)$ that are induced from a homomorphism $\cU_q(\fg)\to \cO_q(\cX)$ to the regular functions of some cluster variety $\cX$, with a suitable polarization of a quantum torus algebra $\cO_q(\cX_\s)$ for some seed $\s$.

It was shown in \cite{Ip18}, \cite{Ip20} and \cite{IM22} that the following list of irreducible regular positive representations is known:
\begin{enumerate}[label=(\arabic*)]
\item The standard representation $\cP_\lambda$,
\item The parabolic representation $\cP_\lambda^J$ for $W_J\subset W$,
\item The degenerate representation $\cP_\lambda^{0,J}$ for $W_J\subset W$ of type $A_{k_1}\x\cdots\x A_{k_m}$,
\item The modular double $\cP_{\wt{\lambda}}^{0,J}$ of $\cP_\lambda^{0,J}$,
\item A mixture of type (2)-(4) for disconnected subsets of Dynkin indices.
\end{enumerate}

By construction, all these representations stem from the standard representation. It is conjectured in \cite[Section 6]{IM22} that the above list is complete.

To give some intuition, recall the isomorphism
\[
\cU_q(\fg)\simeq \cO_q(\sL_{G,\odot})
\]
for simply laced $\fg$ and conjectured for general $\fg$. Based on this isomorphism, it is natural to guess that all irreducible regular positive representations must come from the standard one. However, it is hard to determine if there are other modifications we can do to the standard representation, and generate more irreducible regular positive representations.

Combining Theorem \ref{MainThm}, we conjecture that for any irreducible regular positive representation, the homomorphisms $\cU_q(\fg)\to \cO_q(\cX)$ must send the Chevalley generators to universal polynomials. Even if there are other representations missed by the list, it is interesting to ask whether universal polynomiality always hold.

\medskip

\printbibliography[heading=bibintoc]

\end{document}